\documentclass[11pt,a4paper,twoside]{article}
\usepackage{latexsym,amsmath,amsthm,amssymb,amsfonts}
\usepackage{mathrsfs,multicol}
\usepackage{amscd}
\usepackage{cite}
\usepackage[ansinew]{inputenc}

\def\proof{{\bf Proof.}\quad}
\def\endproof{\hfill\hbox{$\sqcup$}\llap{\hbox{$\sqcap$}}\medskip}
\numberwithin{equation}{section}

\newtheorem{thm}{\indent{\bf Theorem}}[section]

\catcode`@=11

\newskip\plaincentering \plaincentering=0pt plus 1000pt minus 1000pt
\def\@plainlign{\tabskip=0pt\everycr={}}
\def\eqalignno#1{\displ@y \tabskip\plaincentering
  \halign to\displaywidth{\hfil$\@lign\displaystyle{##}$\tabskip\z@skip
    &$\@lign\displaystyle{{}##}$\hfil\tabskip\plaincentering
    &\llap{$\@lign##$}\tabskip\z@skip\crcr
    #1\crcr}}
\def\leqalignno#1{\displ@y \tabskip\plaincentering
  \halign to\displaywidth{\hfil$\@lign\displaystyle{##}$\tabskip\z@skip
    &$\@lign\displaystyle{{}##}$\hfil\tabskip\plaincentering
    &\kern-\displaywidth\rlap{$\@lign##$}\tabskip\displaywidth\crcr
    #1\crcr}}
\def\plainLet@{\relax\iffalse{\fi\let\\=\cr\iffalse}\fi}
\def\plainvspace@{\def\vspace##1{\noalign{\vskip##1}}}

\def\intic@{\mathchoice{\hskip5\p@}{\hskip4\p@}{\hskip4\p@}{\hskip4\p@}}
\def\negintic@
 {\mathchoice{\hskip-5\p@}{\hskip-4\p@}{\hskip-4\p@}{\hskip-4\p@}}
\def\intkern@{\mathchoice{\!\!\!}{\!\!}{\!\!}{\!\!}}
\def\intdots@{\mathchoice{\cdots}{{\cdotp}\mkern1.5mu
    {\cdotp}\mkern1.5mu{\cdotp}}{{\cdotp}\mkern1mu{\cdotp}\mkern1mu
      {\cdotp}}{{\cdotp}\mkern1mu{\cdotp}\mkern1mu{\cdotp}}}
\newcount\intno@
\def\iint{\intno@=\tw@\futurelet\next\ints@}
\def\iiint{\intno@=\thr@@\futurelet\next\ints@}
\def\iiiint{\intno@=4 \futurelet\next\ints@}
\def\idotsint{\intno@=\z@\futurelet\next\ints@}
\def\ints@{\findlimits@\ints@@}
\newif\iflimtoken@
\newif\iflimits@
\def\findlimits@{\limtoken@false\limits@false\ifx\next\limits
 \limtoken@true\limits@true\else\ifx\next\nolimits\limtoken@true\limits@false
    \fi\fi}
\def\multintlimits@{\intop\ifnum\intno@=\z@\intdots@
  \else\intkern@\fi
    \ifnum\intno@>\tw@\intop\intkern@\fi
     \ifnum\intno@>\thr@@\intop\intkern@\fi\intop}
\def\multint@{\int\ifnum\intno@=\z@\intdots@\else\intkern@\fi
   \ifnum\intno@>\tw@\int\intkern@\fi
    \ifnum\intno@>\thr@@\int\intkern@\fi\int}
\def\ints@@{\iflimtoken@\def\ints@@@{\iflimits@
   \negintic@\mathop{\intic@\multintlimits@}\limits\else
    \multint@\nolimits\fi\eat@}\else
     \def\ints@@@{\multint@\nolimits}\fi\ints@@@}
\def\Sb{_\bgroup\vspace@
        \baselineskip=\fontdimen10 \scriptfont\tw@
        \advance\baselineskip by \fontdimen12 \scriptfont\tw@
        \lineskip=\thr@@\fontdimen8 \scriptfont\thr@@
        \lineskiplimit=\thr@@\fontdimen8 \scriptfont\thr@@
        \Let@\vbox\bgroup\halign\bgroup \hfil$\scriptstyle
            {##}$\hfil\cr}
\def\endSb{\crcr\egroup\egroup\egroup}
\def\Sp{^\bgroup\vspace@
        \baselineskip=\fontdimen10 \scriptfont\tw@
        \advance\baselineskip by \fontdimen12 \scriptfont\tw@
        \lineskip=\thr@@\fontdimen8 \scriptfont\thr@@
        \lineskiplimit=\thr@@\fontdimen8 \scriptfont\thr@@
        \Let@\vbox\bgroup\halign\bgroup \hfil$\scriptstyle
            {##}$\hfil\cr}
\def\endSp{\crcr\egroup\egroup\egroup}
\def\Let@{\relax\iffalse{\fi\let\\=\cr\iffalse}\fi}
\def\vspace@{\def\vspace##1{\noalign{\vskip##1 }}}
\def\aligned{\,\vcenter\bgroup\plainvspace@\plainLet@\openup\jot\m@th\ialign
  \bgroup \strut\hfil$\displaystyle{##}$&$\displaystyle{{}##}$\hfil\crcr}
\def\endaligned{\crcr\egroup\egroup}
\def\matrix{\,\vcenter\bgroup\plainLet@\plainvspace@
    \normalbaselines
  \m@th\ialign\bgroup\hfil$##$\hfil&&\quad\hfil$##$\hfil\crcr
    \mathstrut\crcr\noalign{\kern-\baselineskip}}
\def\endmatrix{\crcr\mathstrut\crcr\noalign{\kern-\baselineskip}\egroup
                \egroup\,}
\newtoks\hashtoks@
\hashtoks@={#}
\def\format{\crcr\egroup\iffalse{\fi\ifnum`}=0 \fi\format@}
\def\format@#1\\{\def\preamble@{#1}%
  \def\c{\hfil$\the\hashtoks@$\hfil}%
  \def\r{\hfil$\the\hashtoks@$}%
  \def\l{$\the\hashtoks@$\hfil}%
  \setbox\z@=\hbox{\xdef\Preamble@{\preamble@}}\ifnum`{=0 \fi\iffalse}\fi
   \ialign\bgroup\span\Preamble@\crcr}

\def\cases{\left\{\,\vcenter\bgroup\plainvspace@
     \normalbaselines\openup\jot\m@th
      \plainLet@\ialign\bgroup$\displaystyle{##}$\hfil&
      \quad$\displaystyle{{}##}$\hfil\crcr
      \mathstrut\crcr\noalign{\kern-\baselineskip}}

\newif\iftagsleft@
\tagsleft@true
\def\TagsOnRight{\global\tagsleft@false}
\def\tag#1$${\iftagsleft@\leqno\else\eqno\fi
 \hbox{\def\pagebreak{\global\postdisplaypenalty-\@M}%
 \def\nopagebreak{\global\postdisplaypenalty\@M}\rm(#1\unskip)}%
  $$\postdisplaypenalty\z@\ignorespaces}
\interdisplaylinepenalty=\@M
\def\allowdisplaybreak{\noalign{\allowbreak}}
\def\plainallowdisplaybreak@{\def\allowdisplaybreak{\noalign{\allowbreak}}}
\def\plaindisplaybreak@{\def\displaybreak{\noalign{\break}}}
\def\align#1\endalign{\def\tag{&}\plainvspace@\plainallowdisplaybreak@
\plaindisplaybreak@
  \iftagsleft@\plainlalign@#1\endalign\else
   \plainralign@#1\endalign\fi}
\def\plainralign@#1\endalign{\displ@y\plainLet@\tabskip\plaincentering
\halign to\displaywidth
     {\hfil$\displaystyle{##}$\tabskip=\z@&$\displaystyle{{}##}$\hfil
       \tabskip=\plaincentering&\llap{\hbox{\rm(##\unskip)}}\tabskip\z@\crcr
             #1\crcr}}
\def\plainlalign@
 #1\endalign{\displ@y\plainLet@\tabskip\plaincentering\halign to \displaywidth
   {\hfil$\displaystyle{##}$\tabskip=\z@&$\displaystyle{{}##}$\hfil
   \tabskip=\plaincentering&\kern-\displaywidth
        \rlap{\hbox{\rm(##\unskip)}}\tabskip=\displaywidth\crcr
               #1\crcr}}

\def\re@#1{\par\hangindent\parindent\indent\llap{#1\enspace}\ignorespaces}
\def\qfootnote#1{\edef\@sf{\spacefactor\the\spacefactor}{}#1\@sf
      \insert\footins{\let\egroup=}\footnotesize 
      \interlinepenalty100 \let\par=\endgraf
        \leftskip=0pt \rightskip=0pt
        \splittopskip=10pt plus 1pt minus 1pt \floatingpenalty=20000
   \smallskip\re@{#1}\bgroup\strut\aftergroup{\strut\egroup}\let\next}
\catcode`\@=\active
\TagsOnRight
\begin{document}
\title{\bf Monotonicity of eigenvalues of geometric operaters along the
Ricci-Bourguignon flow\footnote{The research of authors is
supported by NNSFC(no.11471246) and NSFHE(no.KJ2014A257).}}
\author{Fanqi Zeng,\ Qun He,\ Bin Chen \footnote{The corresponding author's
email:~\textsf{Chenbin$@$tongji.edu.cn}\,(B. Chen)}}
\date{}
\maketitle
\begin{quotation}
\noindent{\bf Abstract.}~ In this paper, we study monotonicity of eigenvalues of Laplacian-type operator $-\Delta+cR$, where $c$ is a constant,
 along the Ricci-Bourguignon flow. For $c\neq0$,  We derive monotonicity of the lowest eigenvalue of Laplacian-type operator $-\Delta+cR$ which
generalizes some results of Cao \cite{Cao2007}. For $c=0$,  We derive monotonicity of the first eigenvalue of Laplacian which
generalizes some results of Ma \cite{Ma2006}. Moreover, we prove that when $(M_{3}, g_{0})$ is a closed three manifold with positive Ricci curvature, the eigenvalue of
the Laplacian diverges as $t \rightarrow T$ on a limited maximal
time in terval $[0, T)$, which
generalizes some results of Cerbo and Fabrizio \cite{Fabrizio2007}.
\end{quotation}

\section{Introduction}

Given an $n$-dimensional closed Riemannian manifold $(M,g)$, the
 metric $g=g(t)$ is evolving according to the flow equation
\begin{equation}\label{1Int1}
\frac{\partial}{\partial t}g=-2Ric+2\rho Rg=-2(Ric-\rho Rg),
\end{equation}
with the initial condition $$g(0)=g_{0},$$ where $Ric$ is the Ricci
tensor of the manifold, $R$ is scalar curvature and $\rho$ is a real
constant. When $\rho=\frac{1}{2}$, $\rho=\frac{1}{n}$,
$\rho=\frac{1}{2(n-1)}$ and $\rho=0$, the tensor $Ric-\rho Rg$
corresponds to the Einstein tensor, the traceless Ricci tensor, the
Schouten tensor and the Ricci tensor respectively. Apart these
special values of $\rho$, for which we will call the associated
flows as the name of the corresponding tensor, in general we will
refer to the evolution equation defined by the PDE system
\eqref{1Int1} as the Ricci-Bourguignon flow. Moreover, by a suitable
rescaling in time, when $\rho$ is nonpositive, they can be seen as
an interpolation between the Ricci flow and the Yamabe flow (see
\cite{Brendle2005}, \cite{Ye1994}, for instance), obtained as a
limit when $\rho\rightarrow -\infty$.

The study of these flows was proposed by Jean-Pierre Bourguignon
( see Question 3.24 in \cite{Bourguignon1981}), building on some
unpublished work of Lichnerowicz in the sixties and a paper of Aubin
\cite{Aubin1970}. In 2003, Fischer \cite{Fischer2003} studied a
conformal version of this problem where the scalar curvature is
constrained along the flow. In 2011, Lu, Qing and Zheng
\cite{Lu2011} also proved some results on the conformal
Ricci-Bourguignon flow. Recently, for suitable values of the scalar
parameter involved in these flows, Catino et al. \cite{Catino2015}
proved short time existence and provided curvature estimates and
stated some results on the associated solitons. Precisely, we give some useful conclusions about the Ricci-Bourguignon flow  which will be used later.

\noindent{\bf Proposition 1.1.}  (\cite{Catino2015}) {\it
Under the Ricci-Bourguignon flow equation \eqref{1Int1}, we have
\begin{equation}\label{111Th91}
\frac{\partial}{\partial t}g^{ij}=2(R^{ij}-\rho Rg^{ij}),
\end{equation}
\begin{equation}\label{111Th9}
\frac{\partial}{\partial t}(\,d\upsilon)=(n\rho-1)R\,d\upsilon,
\end{equation}
\begin{equation}\label{1Th11}\aligned
\frac{\partial}{\partial t}(\Gamma^{k}_{ij})
=&-R_{ik,j}-R_{kj,i}+R_{ij,k}\\
&+\rho(\delta^{i}_{k}R_{,j}+\delta^{k}_{j}R_{,i}-g_{ij}R_{,k}),
\endaligned\end{equation}
\begin{equation}\label{111Th92}
\frac{\partial}{\partial t}R=[1-2(n-1)\rho]\Delta R+2|Ric|^{2}-2\rho R^{2}.
\end{equation}}

\noindent{\bf Proposition 1.2.}  (\cite{Catino2015}) {\it ( short time existence)
Let $\rho<\frac{1}{2(n-1)}$. Then, the evolution equation
\eqref{1Int1} has a unique solution for a positive time interval on
any smooth, $n$-dimensional, compact Riemannian manifold $M$
(without boundary) with any initial metric $g_{0}$.}

\noindent{\bf Proposition 1.3.}  (\cite{Catino2015}) {\it ( preserved curvature conditions)
Let $(M,g_{t})_{t\in[0, T)}$ be a compact maximal solution of the
Ricci-Bourguignon flow \eqref{1Int1}. If $\rho\leq\frac{1}{2(n-1)}$,
the minimum of the scalar curvature is nondecreasing along the flow.
In particular, if $R(g_{0})\geq \alpha$, for some $\alpha \in
\mathbb{R}$, then $R(g_{t})\geq \alpha$ for every $t\in[0, T)$.
Morever, if $\alpha >0$ then $T\leq \frac{n}{2(1-n\rho)\alpha}$.}

\noindent{\bf Proposition 1.4.}  (\cite{Catino2015}) {\it ( conditions preserved in $3$-dimensional.)
Let $(M,g_{t})_{t\in[0, T)}$ be a compact, $3$-dimensional, solution
of the Ricci-Bourguignon flow \eqref{1Int1}. If $\rho
 < 1/4$, then \\
 (1)nonnegative Ricci curvature is preserved along the flow;\\
 (2)the pinching inequality $Ric \geq \varepsilon Rg$ is preserved along
the flow for any $\varepsilon \leq \frac{1}{3}$.}

\noindent{\bf Proposition 1.5.}  (\cite{Catino2015}) {\it Let $(M,g_{t})_{t\in[0, T)}$ be a compact solution of the
Ricci-Bourguignon flow \eqref{1Int1} with $\rho\leq\frac{1}{2(n-1)}$,
and such that the initial data $g_{0}$ has nonnegative curvature operator. Then $\widetilde{R}_{g(t)}\geq0$ for every $t\in[0, T)$,
where $\widetilde{R} \in End(\bigwedge^{2}M)$ be the Riemann curvature operator.}

\noindent{\bf Proposition 1.6.}  (\cite{Catino2015}) {\it $\rho<\frac{1}{2(n-1)}$. If $g(t)$ is a compact solution of the Ricci-Bourguignon flow on a maximal
time interval $[0, T)$, $T <+\infty$, then
$$\limsup_{t\rightarrow T}\max_{M}|Riem(\cdot, t)|=+\infty,$$
where $Riem(\cdot, t)$ is Riemann tensor.}

At present, the eigenvalues of geometric operators have become a
powerful tool in the study of geometry and topology of manifolds.
Recently, there has been a lot of work on the eigenvalue problems
under the Ricci flow. On one hand, in \cite{Perelman2002} , Perelman introduced
the so-called $\mathscr{F}$-entropy functional and proved that it is nondecreasing along the Ricci flow coupled to a
backward heat-type equation. The nondecreasing of the functional $\mathscr{F}$ implies the monotonicity of the first
eigenvalue of $-4\Delta+R$ along the Ricci flow.

Cao \cite{Cao2007} extended the operator $$-4\Delta+R,$$ to the following new operator
$$-\Delta+\frac{R}{2}$$ on closed Riemannian manifolds, and showed that the eigenvalues of this new operator are
nondecreasing along the Ricci flow with nonnegative curvature operator. Later Li \cite{Li2007} and
Cao \cite{Cao2008} considered a general operator $$-\Delta+cR,~~~where~~~c\geq\frac{1}{4},$$ and both proved that the first eigenvalue of this operator is nondecreasing along the Ricci flow
without any curvature assumption.
 For the recent
research in this direction, see
\cite{Ling2009,Ling2007,Zhao2013,Guo2013, Zhao2014} and the
references therein.

On the other hand, Ma \cite{Ma2006}
obtained the monotonicity of the first eigenvalue of the Laplacian operator on a domain with Dirichlet
boundary condition along the Ricci flow. Using the
differentiability of the eigenvalues and the corresponding
eigenfunctions of the Laplace operator under the Ricci flow, he
obtained the following result.

\begin{thm}\label{thm1-4}( \cite{Ma2006})
Let $g=g(t)$ be the evolving metric along the Ricci-Hamilton flow
with $g(0)=g_{0}$ being the initial metric in $M$. Let $D$ be a
smooth bounded domain in $(M,g_{0})$. Let $\mu>0$ be the first
eigenvalue of the Laplace operator of the metric $g(t)$. If there is
a constant $a$ such that the scalar curvature $R\geq2a$ in $D \times
\{t\}$ and the Einstein tensor $$E_{ij}\geq-ag_{ij},~~~in ~~~ D
\times \{t\},$$ where $E_{ij}:= R_{ij}-\frac{R}{2}g_{ij}$, then we
have $\frac{d}{dt}\mu\geq 0$, that is, $\mu$ is nondecreasing in
$t$; furthermore, $\frac{d}{dt}\mu(t)>0$ for the scalar curvature
$R$ not being the constant $2a$. The same monotonicity result is
also true for other eigenvalues.
\end{thm}

Motivated by the work of Perelman, Cao and Ma, we shall consider the
eigenvalue of $-\Delta+cR$ with $c$ a constant. For $c\neq0$, inspired by \cite{Cao2007} and \cite{Cao2008},  we can derive the following monotonicity of the lowest eigenvalue of $-\Delta+cR$ under the Ricci-Bourguignon flow equation \eqref{1Int1}. That is, we obtain the

\begin{thm}\label{thm1-11}
Let $(M,g(t))_{t\in[0, T)}$ be a compact maximal solution of the
Ricci-Bourguignon flow \eqref{1Int1}, $g(t)\not\equiv g(0)$ and
$\lambda_{0}(t)$ be
the lowest eigenvalue of the operator $-\Delta+cR$ corresponding to the normalized eigenfunction $f$, that is,
$$(-\Delta+cR)f=\lambda_{0} f,~~~\int_{M} f^2\,d\upsilon=1.$$

(1) If $\rho\leq0$,  $ c\geq\frac{[1-\rho(n-1)]^{2}}{4-8\rho(n-1)}$ and the scalar curvature is nonnegative at the initial time, then the lowest eigenvalue of the operator $-\Delta+cR$ is is strictly monotone increasing in
$[0, T)$ under the Ricci-Bourguignon flow equation \eqref{1Int1}.

(2) If $0<\rho\leq\frac{1}{2(n-1)}$, $c\geq\frac{1-(n-2)\rho}{2[1-2(n-1)\rho]}$ and the curvature operator is nonnegative at the initial time, then the following quantity
\begin{equation}\label{132Th7}
\left(T^{'}-t\right)^{-\alpha}\lambda_{0}(t),
\end{equation}
is strictly monotone increasing
 under the Ricci-Bourguignon flow equation \eqref{1Int1} in $[0, T^{'})$ and $T^{'}=\frac{1}{2(1-\rho)\epsilon}$, where $$\epsilon=\max_{x\in M}R(0),$$ $$\alpha=\frac{2c[1-2(n-1)\rho]+n\rho-1}{2(1-\rho)}\geq0.$$
\end{thm}
\noindent{\bf Remark 1.3.} It should be pointed out that for $\rho = 0$ , our (1) generalizes the corresponding results of Cao in \cite{Cao2007}.

\noindent{\bf Remark 1.4.} It's obvious that (2) will hold whenever the Ricci curvature
is nonnegative, but in general, the nonnegativity of the Ricci curvature is not
preserved. Nevertheless, the nonnegativity of the Ricci curvature is preserved
in dimension three.

\noindent {\bf Corollary 1.1.}  {\it In dimension three, let $g(t)$ and $\lambda_{0}(t)$ be the same as in Theorem 1.2. But here we assume the Ricci curvature  is nonnegative at the initial time.
If $0<\rho\leq\frac{1}{4}$ and $c\geq\frac{1-\rho}{2(1-4\rho)}$, then the following quantity
\begin{equation}\label{132Th7}
\left(T^{'}-t\right)^{-\alpha}\lambda_{0}(t),
\end{equation}
is strictly monotone increasing  under the Ricci-Bourguignon flow equation \eqref{1Int1} in $[0, T^{'})$ and $T^{'}=\frac{1}{2(1-\rho)\epsilon}$, where $$\epsilon=\max_{x\in M}R(0),$$ $$\alpha=\frac{2c(1-4\rho)+3\rho-1}{2(1-\rho)}\geq0.$$ }

For $c= 0$, we derive the following monotonicity of eigenvalues on Laplacian under the Ricci-Bourguignon flow equation \eqref{1Int1}. That is, we obtain the

\begin{thm}\label{thm1-6}
Let $(M,g(t))_{t\in[0, T)}$ be a compact maximal solution of the
Ricci-Bourguignon flow \eqref{1Int1}, $g(t)\not\equiv g(0)$ and $\rho\leq\frac{1}{2(n-1)}$.
Let $\lambda_{1}(t)$ be the first eigenvalue of the Laplace operator of the
metric $g(t)$. If there is a non-negative constant $a$ such that
\begin{equation}\label{1632Th7}
R_{ij}-\frac{1+(2-n)\rho}{2}Rg_{ij}\geq-ag_{ij},~~~in ~~~ M \times [0, T),
\end{equation}
and
\begin{equation}\label{1732Th7}
R\geq\frac{2a}{1-n\rho} , ~~~in ~~~  M \times \{0\},
\end{equation}
then $\lambda_{1}(t)$ is strictly monotone increasing and differentiable almost everywhere along the Ricci-Bourguignon flow in
$[0, T)$.

\end{thm}

\noindent{\bf Remark 1.1.}
(1) In \cite{Wu2011}, Wu et al. proved a similar result about the $p$-Laplace operator along the Ricci flow, where they
assumed $R \geq a p$ and $R \not\equiv a p$ in $M \times \{0\}$, which are a little stronger than assumptions of Theorem 1.3. The key difference is that we use Lemma 2.3.

(2) It should be pointed out that for $\rho = 0$ , the above theorem is similar to Ma's main result for the first eigenvalue of
the Laplace operator in \cite{Ma2006}. Moreover, our assumptions are weaker than Ma's.

(3) If $a<0$, there doesn't exist any scalar curvature which satisfies \eqref{1632Th7} and \eqref{1732Th7}.

Cerbo and Fabrizio in \cite{Fabrizio2007} proved that when $(M_{3},
g_{0})$ is a closed manifold with positive Ricci curvature,
the eigenvalues of the Laplacian along the Ricci flow diverges as $t
\rightarrow T$. Motivated by the work of Cerbo and Fabrizio,  we can obtain a similar result under the
Ricci-Bourguignon flow. That is, we obtain the

\begin{thm}\label{thm1-13}
Let $\lambda$ be the eigenvalue of the Laplacian along the Ricci-Bourguignon, that is, $-\Delta f=\lambda f$.
Let $(M_{3}, g(t))$ be a compact maximal solution of the Ricci-Bourguignon flow on a
closed $3$-manifold with positive Ricci curvature initially and
$\rho<\frac{1}{4}$, then
$$\lim\limits_{t\rightarrow T}\lambda(t)=\infty$$
on a limited maximal
time in terval $[0, T)$.
\end{thm}

\noindent{\bf Remark 1.2.} It should be pointed out that for $\rho = 0$ , our above
results on eigenvalues reduce to the Proposition 3.1 of Cerbo and Fabrizio in \cite{Fabrizio2007}.


\section{Preliminaries}

In this section, we will first give
the definition for the first eigenvalue (the lowest eigenvalue) of the Laplace operator (Laplacian-type operator $-\Delta+cR$) under the Ricci-Bourguignon flow on a closed manifold. Then, we will show that the first eigenvalue (the lowest eigenvalue) of the Laplace operator (Laplacian-type operator $-\Delta+cR$) is a
continuous function along the Ricci-Bourguignon flow. Finally, under the Ricci-Bourguignon flow, we show that if $R(g_{0})\geq \beta$, for some $\beta \in
\mathbb{R}$, and $g(t)\not\equiv g(0)$, then $\max R(g_{t})> \beta$ for every $t\in[0, T)$.

Throughout, $M$ will be taken to be a closed manifold (i.e., compact without boundary). We use moving frames in all calculations and adopt the
following index convention:
$$1\leq i,j,k,\cdots\leq n$$
throughout this paper.

Now we recall the definition of the first eigenvalue of the Laplace
operator on a closed manifold under the Ricci-Bourguignon flow. Let
$(M_{n}, g(t))$ be a solution of the Ricci-Bourguignon flow on the
time interval $[0, T)$. Consider the first nonzero eigenvalue of the
Laplace operator at time $t$, where $0\leq t<T$,
$$\lambda_{1}(t)=\inf\left\{\int_{M}|\nabla f|^2\,d\upsilon: f\in
W^{1,2},\int_{M}f^2\,d\upsilon=1 ~~~and~~~ \int_{M} f
\,d\upsilon=0\right\},$$ where $d\upsilon$ denotes the volume form
of the metric $g=g(t)$. Meanwhile the corresponding eigenfunction
$f$ satisfies the equation
$$-\Delta f(t)=\lambda_{1}(t) f(t),$$
where $\Delta$ is the Laplace operator with respect to $g(t)$, given
by
$$\Delta_{g(t)}=\frac{1}{\sqrt{|g(t)|}}\partial_{i}
(\sqrt{|g(t)|}g(t)^{ij}\partial_{j}),$$ and
$g(t)^{ij}=g(t)_{ij}^{-1}$ is the inverse of the matrix $g(t)$ and
$|g|=det(g_{i j})$.

Note that it is not clear whether the first eigenvalue or the
corresponding eigenfunction of the Laplace operator is
$C^{1}$-differentiable under the Ricci-Bourguignon flow. When $\rho=
0$, where the Ricci-Bourguignon flow is the Ricci-Hamilton flow,
many papers have pointed out that its differentiability under the
Ricci-Hamilton flow follows from eigenvalue perturbation theory
(e.g., see \cite{Kato1984,Kleiner2008,Reed1978}). But $\rho\neq0$,
as far as we are aware, the differentiability of the first
eigenvalue and eigenfuntion of the Laplace operator under the
Ricci-Bourguignon flow has not been known until now. So we can not
use Ma's trick to derive the monotonicity of the first eigenvalue of
the Laplace operator. Although, we do not know the differentiability for $\lambda_{1}(t)$, following the techniques of Wu et al. \cite{Wu2011}, we will see that $\lambda_{1}(t)$ in fact is a continuous function along the Ricci-Bourguignon flow on [0, T ).

\noindent{\bf Lemma 2.1.}(\cite{Wu2011}) {\it If $g_{1}$ and $g_{2}$ are two metrics which satisfy
$$(1+\varepsilon)^{-1}g_{1}\leq g_{2}\leq(1+\varepsilon)g_{1},$$
then, we have
\begin{equation}\label{1001Th91}
(1+\varepsilon)^{-(n+1)}\leq\frac{\lambda_{1}(g_{1})}{\lambda_{1}(g_{2})}\leq(1+\varepsilon)^{(n+1)}.
\end{equation}
In particular, $\lambda_{1}((g(t)))$ is a continuous function in the $t$-variable.}

Next we recall the definition of the lowest eigenvalue of $-\Delta+cR$. Let $\lambda_{0}(t)$ be the lowest eigenvalue of $-\Delta+cR$. Given a metric $g$ on a closed manifold $M^{n}$, we define the functional  $ \lambda_{0}$ by
\begin{equation}\label{1111Th91}
\lambda_{0}(t)=\inf\left\{\mathscr{G}(g, f): \int_{M}f^2\,d\upsilon=1, f>0,~~~ f\in
W^{1,2}\right\},
\end{equation}
where $\mathscr{G}(g, f)=\int_{M}(f(-\Delta f)+cRf^{2})\,d\upsilon$.

We also do not know the differentiability for $\lambda_{0}(t)$. But,
following the techniques of \cite{Chow2008}, we will see that $\lambda_{0}(t)$ in fact is a continuous function along the Ricci-Bourguignon flow on $[0, T )$.

\noindent{\bf Lemma 2.2.}  (\cite{Chow2008}) {\it If $g_{1}$ and $g_{2}$ are two metrics on $M$ which satisfy
$$(1+\varepsilon)^{-1}g_{1}\leq g_{2}\leq(1+\varepsilon)g_{1} ~~~and~~~R(g_{1})-\varepsilon\leq R(g_{2})\leq R(g_{1})+\varepsilon, $$
then
\begin{equation}\label{zTh91}\aligned
 &\lambda_{0}(g_{2})-\lambda_{0}(g_{1})\\
 \leq&\left((1+\varepsilon)^{\frac{n}{2}+1}-(1+\varepsilon)^{-\frac{n}{2}}\right)(1+\varepsilon)^{\frac{n}{2}}(\lambda_{0}(g_{1})-|c|\min R(g_{1}) )\\
&+\left((1+\delta)|c|\max|R(g_{2})-R(g_{1})|+2\delta|c|\max|R(g_{1})|\right)(1+\varepsilon)^{\frac{n}{2}},
 \endaligned\end{equation}
where $\delta\rightarrow 0$ as $\varepsilon\rightarrow 0$. In particular, $\lambda_{0}$ is a continuous function with respect to the $C^{2}$-topology.}

At last, we present the following lemma.

\noindent{\bf Lemma 2.3.} {\it Let $(M,g_{t})_{t\in[0, T)}$ be a compact maximal solution of the
Ricci-Bourguignon flow \eqref{1Int1}. If $\rho\leq\frac{1}{2(n-1)}$,
 $R(0)\geq \beta$, for some $\beta \in
\mathbb{R}$, and $g(t)\not\equiv g(0)$, then $\max R(t)> \beta$ for every $t\in(0, T)$.}

\proof  From Proposition 1.3, we know that $R(t)\geq \beta$ for every $t\in[0, T)$. If $\max R(t_{0})=\beta$ for some $t_{0}\in(0, T)$, we have $R(t_{0})\equiv\beta$. Since \eqref{111Th92} and $\frac{\partial R}{\partial t}\leq0$, we have $$\frac{1}{n}R^{2}(t_{0})\leq|Ric|^{2}(t_{0})\leq\rho R^{2}(t_{0})\leq\frac{1}{2(n-1)}R^{2}(t_{0}).$$
Obviously, we have $R(t_{0})=0$ and $Ric(t_{0})=0$. Hence, $\max R(t_{0})=\beta=0$. Therefore, if $\beta\neq0$, we have $\max R(t)> \beta$ for every $t\in[0, T)$.

When $\beta=0$, let $I=\{t>0: \max R(t)>0\}$. If $I=\emptyset$, then we have $R(t)\equiv0$ and $Ric(t)\equiv0$. Hence, we have $g(t)= g(0)$ which is in contradiction with $g(t)\not\equiv g(0)$. When $I\neq\emptyset$ and $t_{1}\in I$, for any $t_{0}$ and $0<t_{0}<t_{1}$, if $\max R(t_{0})=0$, then $R(t_{0})\equiv0$ and $Ric(t_{0})\equiv0$. Hence, in $[t_{0}, T)$, $g(t_{1})= g(t_{0})$. So we have $Ric(t_{1})= Ric(t_{0})=0$ which is in contradiction with $\max R(t_{1})>0$. Hence, $t_{0}\in I$. Since $t_{0} \in(0, t_{1})$ is arbitrary, we have $(0, t_{1}]\subset I$. By the strong maximum principle, we have $(0, T)\subset I$. Therefore we finish the proof of Lemma 2.3.

\section{Proof of Theorems 1.2 }

In this section, we will prove Theorem 1.2 in Sect.1. In order to achieve this, we first prove
the following two lemmas. Our proof uses some tricks in  \cite{Cao2007} and \cite{Cao2008}.

Let $M$ be an $n$-dimensional closed Riemannian manifold, and $g(t)$
be a smooth solution of the Ricci-Bourguignon flow on the time
interval $[0, T)$. Let $\lambda_{0}(t)$ be
the lowest eigenvalue of the operator $-\Delta+cR$ corresponding to the normalized eigenfunction $f$, that is,
$$(-\Delta+cR)f=\lambda_{0} f,~~~\int_{M} f^2\,d\upsilon=1.$$
For any $t_{0}\in[0, T)$,  there existed a smooth function $\varphi(t)>0$
satisfying
\begin{equation}\label{781Th2}
\int_{M}\varphi^2(t)\,d\upsilon=1,
\end{equation}
 and $\varphi(t_{0})=f(t_{0})$.
Let
\begin{equation}\label{71Th2}
\mu(t)=\int_{M}(\varphi(t)(-\Delta \varphi(t))+cR\varphi^{2}(t))\,d\upsilon,
\end{equation}
then $\mu(t)$ is a smooth function by definition. And at time $t_{0}$, we conclude that $$\mu(t_{0})=\lambda_{0}(t_{0}).$$

We first give the following Lemma.

\noindent{\bf Lemma 3.1.} {\it
Let $g(t)$, $t\in[0,T)$, be a
solution of the Ricci-Bourguignon flow on an n-dimensional closed
manifold $M$, and let $\lambda_{0}(t)$ be the lowest eigenvalue of $-\Delta+cR$ under the Ricci-Bourguignon flow.
Assume that
$f(t_{0})$ is the corresponding eigenfunction of $\lambda_{0}(t)$ at time
$t_{0}\in[0,T)$.
Let $\mu(t)$ be a smooth function defined by \eqref{71Th2}. Then
we have
\begin{equation}\label{1Th109}\aligned
\frac{d}{dt}\mu(t)\mid_{t=t_{0}}=&(A-2\rho)c\int_{M}R^{2} f^{2}\,d\upsilon+(A-2\rho)\int_{M}R |\nabla f|^{2}\,d\upsilon\\
 &-A\lambda_{0}\int_{M}Rf^{2}\,d\upsilon
 +2\int_{M}Ric(\nabla f,\nabla f)d\upsilon\\
 &+2c\int_{M}|Ric|^{2}f^{2}d\upsilon,
  \endaligned\end{equation}
 where $A=-1+n\rho+2c[1-2(n-1)\rho]$.
}

\proof  The proof is by straightforward computation.
Notice that
\begin{equation}\label{1Th9}
\frac{\partial}{\partial t}(\Delta \phi)
=2R_{ij}\phi_{ij}+\Delta (\phi_{t})-2\rho R\Delta \phi-(2-n)\rho
R_{,k}\phi_{k}.
\end{equation}
Using Proposition 1.1, we have
\begin{equation}\label{72Th109}\aligned
\frac{d}{dt}\mu(t)\mid_{t=t_{0}}=&\int_{M}\partial_{t}(-\Delta \varphi+cR\varphi)\varphi\,d\upsilon+\int_{M}(-\Delta \varphi+cR\varphi)\partial_{t}(\varphi\,d\upsilon)\\
 =&\int_{M}[\partial_{t}(-\Delta \varphi)+c\varphi\partial_{t}R+cR\partial_{t}\varphi]\varphi\,d\upsilon
 +\int_{M}(-\Delta \varphi+cR\varphi)\partial_{t}(\varphi\,d\upsilon)\\
  =&\int_{M}[-2R_{ij}\varphi_{ij}-\Delta (\varphi_{t})+2\rho R\Delta \varphi+(2-n)\rho
R_{,k}\varphi_{k}]\varphi\,d\upsilon\\
&+\int_{M}[c\varphi\partial_{t}R+cR\partial_{t}\varphi]\varphi\,d\upsilon+\int_{M}(-\Delta \varphi+cR\varphi)\partial_{t}(\varphi\,d\upsilon)\\
=&\int_{M}[-2R_{ij}\varphi_{ij}\varphi+2\rho R\varphi\Delta \varphi+(2-n)\rho R_{,k}\varphi_{k}\varphi]\,d\upsilon\\
&+c\int_{M}\{[1-2(n-1)\rho]\Delta R+2|Ric|^{2}-2\rho R^{2}\}\varphi^{2}\,d\upsilon\\
&+\int_{M}(-\Delta \varphi+cR\varphi)[\partial_{t}(\varphi)d\upsilon+\partial_{t}(\varphi\,d\upsilon)].
  \endaligned\end{equation}
From $R_{,i}=2R_{ij,j}$ and the Stokes formula, we have
\begin{equation}\label{73Th2}
\int_{M}\varphi^{2}\Delta R\,d\upsilon=\int_{M}2 R(|\nabla \varphi|^{2}+\varphi\Delta \varphi)\,d\upsilon,
\end{equation}
\begin{equation}\label{74Th2}
\int_{M}R_{,k}\phi_{k}\phi\,d\upsilon=\int_{M}2 R(|\nabla \phi|^{2}+\phi\Delta \phi)\,d\upsilon,
\end{equation}
and
\begin{equation}\label{111Th91}\aligned
 \int_{M} -R_{ij}\phi_{ij}\phi\,d\upsilon=&\int_{M} (R_{ij}\phi)_{j}\phi_{i}\,d\upsilon\\
 =&\int_{M} R_{ij,j}\phi\phi_{i}\,d\upsilon+\int_{M} R_{ij}\phi_{j}\phi_{i}\,d\upsilon\\
 =&\frac{1}{2}\int_{M} R_{,i}\phi\phi_{i}\,d\upsilon+\int_{M} R_{ij}\phi_{j}\phi_{i}\,d\upsilon\\
=&-\frac{1}{2}\int_{M} R(\phi\phi_{i})_{i}\,d\upsilon+\int\limits R_{ij}\phi_{j}\phi_{i}\,d\upsilon\\
=&-\frac{1}{2}\int_{M} R\Delta \phi\phi\,d\upsilon-\frac{1}{2}\int_{M}
R|\nabla \phi|^{2}\,d\upsilon+\int_{M} R_{ij}\phi_{j}\phi_{i}\,d\upsilon.
\endaligned\end{equation}

On the other hand, at time $t_{0}$, $\varphi$ is the eigenfunction of $\lambda_{0}(t_{0})$, i.e., $(-\Delta+cR)\varphi=\lambda_{0}\varphi$, we have
\begin{equation}\label{76Th2}
\int_{M}(-\Delta \varphi+cR\varphi)[\partial_{t}(\varphi)d\upsilon+\partial_{t}(\varphi\,d\upsilon)]=\lambda_{0}\int_{M} \varphi[\partial_{t}(\varphi)d\upsilon+\partial_{t}(\varphi\,d\upsilon)]=0.
\end{equation}
The last equality holds because of \eqref{781Th2}.
Inserting \eqref{73Th2}, \eqref{74Th2}, \eqref{111Th91} and \eqref{76Th2} into \eqref{72Th109}, at $t=t_{0}$, yields
\begin{equation}\label{1Th189}\aligned
\frac{d}{dt}\mu(t)\mid_{t=t_{0}}=&\left(-1+n\rho+2c[1-2(n-1)\rho]\right)\int_{M}R\varphi\Delta \varphi\,d\upsilon\\
 &+\left(-1+(n-2)\rho+2c[1-2(n-1)\rho]\right)\int_{M}R |\nabla \varphi|^{2}\,d\upsilon\\
 &+2\int_{M}R_{ij}\varphi_{i}\varphi_{j}d\upsilon+2c\int_{M}|Ric|^{2}\varphi^{2}d\upsilon-2c\rho\int_{M}R\varphi^{2}\,d\upsilon.
  \endaligned\end{equation}
Inserting  $\Delta \varphi=cR\varphi-\lambda_{0} \varphi$ into \eqref{1Th189}, at $t=t_{0}$, so that
\begin{equation}\label{1Th1809}\aligned
 \frac{d}{dt}\mu(t)\mid_{t=t_{0}}=&(A-2\rho)c\int_{M}R^{2} \varphi^{2}\,d\upsilon+(A-2\rho)\int_{M}R |\nabla \varphi|^{2}\,d\upsilon\\
 &-A\mu\int_{M}R\varphi^{2}\,d\upsilon
 +2\int_{M}R_{ij}\varphi_{i}\varphi_{j}d\upsilon\\
 &+2c\int_{M}|Ric|^{2}\varphi^{2}d\upsilon.
  \endaligned\end{equation}
Therefore we finish the proof of Lemma 3.1.

Then we give the second Lemma.

\noindent{\bf Lemma 3.2.} {\it
Let $g(t)$, $t\in[0,T)$, be a
solution of the Ricci-Bourguignon flow on an n-dimensional closed
manifold $M$, and let $\lambda_{0}(t)$ be the lowest eigenvalue of $-\Delta+cR$ under the Ricci-Bourguignon flow.
Assume that
$f(t_{0})$ is the corresponding eigenfunction of $\lambda_{0}(t)$ at time
$t_{0}\in[0,T)$.
Let $\mu(t)$ be a smooth function defined by \eqref{71Th2}. Then
we have
\begin{equation}\label{1Th1009}\aligned
\frac{d}{dt}\mu(t)\mid_{t=t_{0}}=&\frac{\left(1-\rho(n-1)\right)^{2}}{2-4\rho(n-1)}\int_{M}|R_{ij}-2\frac{1-2\rho(n-1)}{1-\rho(n-1)}(\log f)_{ij}|^{2}f^{2}\,d\upsilon\\
&+\left(2c-\frac{\left(1-\rho(n-1)\right)^{2}}{2-4\rho(n-1)}\right)\int_{M} |Ric|^{2}f^{2}\,d\upsilon-\rho\lambda_{0}\int_{M} Rf^{2}\,d\upsilon\\
 &-\rho c\int_{M}R^{2}f^{2}d\upsilon
 -\rho\int_{M}R|\nabla f|^{2}d\upsilon.
  \endaligned\end{equation}
}

\proof  The proof is by straightforward computation.
\begin{equation}\label{h1Th109}\aligned
&\int_{M}|R_{ij}-2k(\log f)_{ij}|^{2}f^{2}\,d\upsilon\\
=&\int_{M} |Ric|^{2}f^{2}\,d\upsilon+4k^{2}\int_{M} |\nabla^{2}(\log f)|^{2}f^{2}\,d\upsilon-4k\int_{M} R_{ij}(\log f)_{ij}f^{2}\,d\upsilon,
  \endaligned\end{equation}
where $k=\frac{1-2\rho(n-1)}{1-\rho(n-1)}$.
 From \cite{Cao2008}, we can get
\begin{equation}\label{h2Th109}\aligned
&4k^{2}\int_{M} |\nabla^{2}(\log f)|^{2}f^{2}\,d\upsilon\\
=&2k^{2}c\int_{M} R\Delta f^{2}\,d\upsilon-4k^{2}\int_{M} R_{ij}f_{i}f_{j}\,d\upsilon\\
=&4k^{2}c\int_{M} R(f\Delta f+|\nabla f|^{2})\,d\upsilon-4k^{2}\int_{M} R_{ij}f_{i}f_{j}\,d\upsilon
  \endaligned\end{equation}
  and
\begin{equation}\label{h3Th109}\aligned
&-4k\int_{M} R_{ij}(\log f)_{ij}f^{2}\,d\upsilon\\
=&-2k\int_{M} R(f\Delta f+|\nabla f|^{2})\,d\upsilon+8k\int_{M} R_{ij}f_{i}f_{j}\,d\upsilon.
  \endaligned\end{equation}
  Combining \eqref{h2Th109} and \eqref{h3Th109}, we arrive at
\begin{equation}\label{h4Th109}\aligned
&\int_{M}|R_{ij}-2k(\log f)_{ij}|^{2}f^{2}\,d\upsilon\\
=&\int_{M} |Ric|^{2}f^{2}\,d\upsilon+(8k-4k^{2})\int_{M} R_{ij}f_{i}f_{j}\,d\upsilon\\
&+2k(2kc-1)\left(c\int_{M}R^{2}f^{2}d\upsilon-\lambda_{0}\int_{M}Rf^{2}d\upsilon\right)\\
&+2k(2kc-1)\int_{M}R|\nabla f|^{2}d\upsilon.
  \endaligned\end{equation}
Multipling $\frac{1}{2k(2-k)}$ on the both sides of \eqref{h4Th109} , we
conclude that
\begin{equation}\label{h5Th109}\aligned
&\frac{1}{2k(2-k)}\int_{M}|R_{ij}-2k(\log f)_{ij}|^{2}f^{2}\,d\upsilon\\
=&\frac{1}{2k(2-k)}\int_{M} |Ric|^{2}f^{2}\,d\upsilon+2\int_{M} R_{ij}f_{i}f_{j}\,d\upsilon\\
&+\frac{2kc-1}{2-k}\left(c\int_{M}R^{2}f^{2}d\upsilon-\lambda_{0}\int_{M}Rf^{2}d\upsilon\right)\\
&+\frac{2kc-1}{2-k}\int_{M}R|\nabla f|^{2}d\upsilon.
  \endaligned\end{equation}
Let \eqref{1Th109} subtract \eqref{h5Th109}, we have
\begin{equation}\label{h6Th109}\aligned
\frac{d}{dt}\mu(t)\mid_{t=t_{0}}=&\frac{1}{2k(2-k)}\int_{M}|R_{ij}-2k(\log f)_{ij}|^{2}f^{2}\,d\upsilon\\
&+\left(2c-\frac{1}{2k(2-k)}\right)\int_{M} |Ric|^{2}f^{2}\,d\upsilon\\
 &+\left(-A+\frac{2kc-1}{2-k}\right)\lambda_{0}\int_{M} Rf^{2}\,d\upsilon\\
 &+\left(A-2\rho-\frac{2kc-1}{2-k}\right)c\int_{M}R^{2}f^{2}d\upsilon\\
 &+\left(A-2\rho-\frac{2kc-1}{2-k}\right)\int_{M}R|\nabla f|^{2}d\upsilon.
  \endaligned\end{equation}
Inserting  $A=-1+n\rho+2c[1-2(n-1)\rho]$ and  $k=\frac{1-2\rho(n-1)}{1-\rho(n-1)}$ into \eqref{h6Th109}, we obtain  \eqref{1Th109}.
Therefore we finish the proof of Lemma 3.2.

\noindent{\bf Proof of Theorems 1.2.}   We first proof (1). If $\rho\leq0$,  $ c\geq\frac{[1-\rho(n-1)]^{2}}{4-8\rho(n-1)}$ and $R\geq0$ in $M\times \{0\}$, from Proposition 1.3 and Lemma 2.3, we know that $\max R(g_{t})> 0$ for every $t\in(0, T)$ and
$$\frac{\left(1-\rho(n-1)\right)^{2}}{2-4\rho(n-1)}>0.$$
 By \eqref{1Th1009},  we obtain
 \begin{equation}\label{1Th83}
\frac{d}{dt}\mu(t)\mid_{t=t_{0}}\geq0.
\end{equation}
Since the eigenfunction of the lowest eigenvalue is not equal to $0$ almost everywhere along the Ricci-Bourguignon flow and $\mu(t)$ is a smooth function with respect to $t$-variable,
we have
\begin{equation}\label{1Th84}
\frac{d}{dt}\mu(t)>0,
\end{equation}
in $(t_{0}-\delta, t_{0}+\delta)$, where $\delta>0$ is sufficiently small.
So we get
\begin{equation}\label{1Th85}
\mu( t_{0})>\mu( t_{1})
\end{equation}
for any  $ t_{1}\in (t_{0}-\delta, t_{0}+\delta)$  and $ t_{1}< t_{0}$.

Notice that $$\mu(t_{0})=\lambda_{0}(t_{0})~~~~~ and ~~~~~\mu(t_{0})\geq\lambda_{0}(t_{1}).$$ This implies $\lambda_{0}(t_{0})>\lambda_{0}(t_{1})$ for any $t_{0}> t_{1}$.
Since $t_{0} \in [0, T)$ is arbitrary, $\lambda_{0}(t)$ is strictly monotone increasing in
$[0, T)$.
Therefore we finish the proof of (1).

Next  we proof (2). By Proposition 1.5,  we know that the nonnegativity of the
curvature operator is preserved by the Ricci-Bourguignon flow. This implies that the Ricci
curvature is also nonnegative, and we have $|Ric|^{2}\leq R^{2}$.
  The evolution equation  of scalar curvature satisfies
\begin{equation}\label{181Th92}\aligned
\frac{\partial}{\partial t}R=&[1-2(n-1)\rho]\Delta R+2|Ric|^{2}-2\rho R^{2}\\
&\leq[1-2(n-1)\rho]\Delta R+2(1-\rho) R^{2}.
\endaligned\end{equation}
Let $\sigma(t)$ be the solution of the following ODE with initial value
\begin{equation}\label{Zeng8}
\left\{\begin{array}{l}
\frac{\partial\sigma(t)}{\partial t}=2(1-\rho)\sigma^{2}, \\
\sigma(0)=\max_{x\in M}R(0).
\end{array}\right.
\end{equation}
By the maximum principle, let $\epsilon=\max_{x\in M}R(0)$, we can get
$$R(t)\leq\sigma(t)=\left(-2(1-\rho)t+\frac{1}{\epsilon}\right)^{-1}$$ on $[0, T^{''})$, where $T^{''}=\min\{T^{'}, T\}$ and $T^{'}=\frac{1}{2(1-\rho)\epsilon}$.  Because of Lemma 2.6, we have $T^{'}\leq T$. Hence $R(t)\leq\sigma(t)$ on $[0, T^{'})$.
If $0<\rho\leq\frac{1}{2(n-1)}$ and $c\geq\frac{1-(n-2)\rho}{2[1-2(n-1)\rho]}$, we have $(A-2\rho)c\geq0$ and $A\geq0$. Since the eigenfunction of the first eigenvalue is not equal to $0$ almost everywhere along the Ricci-Bourguignon flow, from Lemma 2.3 and \eqref{1Th1809}, we have
\begin{equation}\label{Zeng3}\aligned
\frac{d}{dt}\mu\mid_{t=t_{0}}\geq&-A\mu\int_{M}Rf^{2}\,d\upsilon\\
>&-A\mu\sigma,
\endaligned\end{equation}
which implies
$$(\frac{d}{dt}\mu-A\mu\sigma)\mid_{t=t_{0}}>0.$$

Using the similar arguments as (1), we know that $$\left(T^{'}-t\right)^{-\alpha}\lambda_{0}(t)$$ is strictly monotone increasing under the Ricci-Bourguignon flow equation \eqref{1Int1} on $[0, T^{'})$ and $T^{'}=\frac{1}{2(1-\rho)\epsilon}$, where $$\epsilon=\max_{x\in M}R(0),$$ $$\alpha=\frac{2c[1-2(n-1)\rho]+n\rho-1}{2(1-\rho)}\geq0,$$
 which shows (2) holds. Therefore we finish the proof of Theorem 1.2.

\section{Proof of Theorems 1.3 }

In this section, we will prove Theorem 1.3 in Sect.1. In order to achieve this, we first prove
the following lemma. Our proof involves choosing a proper smooth function, which
seems to be a delicate trick.

 Let $M$ be an $n$-dimensional closed Riemannian manifold, and $g(t)$
be a smooth solution of the Ricci-Bourguignon flow on the time
interval $[0, T)$. Let $\lambda_{1}(t)$ be the first eigenvalue of the
Laplace operator under the Ricci-Bourguignon flow and
$f(t_{0})$ be the corresponding eigenfunction of $\lambda_{1}(t)$ at time
$t_{0}\in[0,T)$, i.e.,
\begin{equation}\label{1Th7}
-\Delta_{g(t_{0})} f(t_{0})=\lambda_{1}(t_{0}) f(t_{0}).
\end{equation} For any $t_{0}\in[0, T)$,  Wu in \cite{Wu2011} pointed out that there existed a smooth function
$$ \phi(t)=\frac{\psi(t)}{(\int_{M} \psi(t)^2\,d\upsilon)^{\frac{1}{2}}},~~~where~~~~~\psi(t)=f(t_{0})\left(\frac{det(g_{ij}(t_{0}))}{det(g_{ij}(t))}\right)^{\frac{1}{2}}$$
satisfying
\begin{equation}\label{11Th7}\int_{M} \phi(t)^2\,d\upsilon=1~~~~~ and ~~~~~\int_{M}\phi(t)\,d\upsilon=0,\end{equation}
and $\varphi(t_{0})=f(t_{0})$.
Now we define a general smooth function as follows:
\begin{equation}\label{12Th7}\mu(t)=\int_{M}\phi(t)(-\Delta \phi(t))\,d\upsilon.\end{equation}
In general, $\mu(t)$ is not equal to $\lambda_{1}(t)$. But at time $t_{0}$, we conclude that $$\mu(t_{0})=\lambda_{1}(t_{0}).$$

\noindent{\bf Lemma 4.1.} {\it Let $g(t)$, $t\in[0,T)$, be a
solution of the Ricci-Bourguignon flow on an n-dimensional closed
manifold $M$, and let $\lambda_{1}(t)$ be the first eigenvalue of the
Laplace operator under the Ricci-Bourguignon flow. Assume that
$f(t_{0})$ is the corresponding eigenfunction of $\lambda_{1}(t)$ at time
$t_{0}\in[0,T)$, i.e.,
\begin{equation}\label{1Th7}
-\Delta_{g(t_{0})} f(t_{0})=\lambda_{1}(t_{0}) f(t_{0}).
\end{equation}
Let $\mu(t)$ be a smooth function defined by \eqref{12Th7}. Then
we have
\begin{equation}\label{1Th8}\aligned
\frac{d}{dt}\mu(t)\mid_{t=t_{0}}=\int_{M}
\{2R_{ij}f_{i}f_{j}&+(1-n\rho)\lambda_{1}Rf^{2}\\
&-[(2-n)\rho+1]R|\nabla
f|^{2}\}\,d\upsilon.
\endaligned\end{equation}
}

\proof The proof is by direct computations.
\begin{equation}\label{1Th19}\aligned
 \frac{d}{dt}\mu(t)\mid_{t=t_{0}}=&\int_{M}\partial_{t}(-\phi\Delta \phi)\,d\upsilon+\int_{M}(-\phi\Delta \phi)\partial_{t}(\,d\upsilon)\\
 =&\int_{M} [-2R_{ij}\phi_{ij}-\Delta (\partial_{t}\phi)+2\rho R\Delta \phi+(2-n)\rho
R_{,k}\phi_{k}]\phi\,d\upsilon\\
&+\int_{M}(-\Delta \phi)\partial_{t}\phi\,d\upsilon
+\int_{M}(-\Delta \phi)\phi(n\rho-1)R\,d\upsilon\\
=&\int_{M} -2R_{ij}\phi_{ij}\phi\,d\upsilon+\int_{M}-2(\Delta
\phi)\partial_{t}\phi\,d\upsilon\\
&+(2-n)\rho\int_{M} R_{,k}\phi_{k}\phi\,d\upsilon+
[1+(2-n)\rho]\int_{M}R(\Delta \phi)\phi\,d\upsilon.
\endaligned\end{equation}

From \eqref{74Th2} and \eqref{111Th91}, we have
\begin{equation}\label{71Th19}\aligned
\frac{d}{dt}\mu(t)\mid_{t=t_{0}}=&-\int_{M} R\Delta \phi\phi\,d\upsilon-\int_{M} R|\nabla
\phi|^{2}\,d\upsilon+2\int_{M}
R_{ij}\phi_{j}\phi_{i}\,d\upsilon\\
&+\int_{M}-2(\Delta
\phi)\partial_{t}\phi\,d\upsilon
-(2-n)\rho\int_{M} R\Delta \phi\phi\,d\upsilon\\
&-(2-n)\rho\int_{M} R|\nabla
\phi|^{2}\,d\upsilon
+[1+(2-n)\rho]\int_{M}R(\Delta \phi)\phi\,d\upsilon\\
=&2\int_{M} R_{ij}\phi_{j}\phi_{i}\,d\upsilon+2\mu\int_{M}
\phi\partial_{t}\phi\,d\upsilon\\
&-[1+(2-n)\rho]\int_{M}R|\nabla \phi|^{2}\,d\upsilon.
\endaligned\end{equation}

Under the Ricci-Bourguignon flow, from the constraint condition
\eqref{11Th7}, we can get
\begin{equation}\label{1Th119}
2\int_{M} \phi\partial_{t}\phi d\upsilon=-(n\rho-1)\int_{M} \phi^{2}Rd\upsilon.
\end{equation}
Hence, at time $t_{0}$, the desired proposition follows from
substituting \eqref{1Th119} into the \eqref{71Th19}. Therefore we finish the proof of Lemma 3.1.
\endproof

\noindent{\bf Proof of Theorems 1.3.}  We assume that for any time $t_{0} \in [0, T)$, if $f(t_{0})$ is the corresponding
eigenfunction of the first eigenvalue $\lambda_{1}(t_{0})$, then we have $\lambda_{1}(t_{0})= \mu( t_{0})$.
By Lemma 4.1, we have
\begin{equation}\label{1Th81}\aligned
\frac{d}{dt}\mu(t)\mid_{t=t_{0}}=&\int_{M} \{(1-n\rho)\lambda_{1}
Rf^{2}+2R_{ij}f_{i}f_{j}\\
&-[(2-n)\rho+1]R|\nabla f|^{2}\}\,d\upsilon\\
=&\int_{M}\{2R_{ij}
-[(2-n)\rho+1]Rg_{ij}\}f_{i}f_{j}\,d\upsilon\\
&+\int_{M} (1-n\rho)\lambda_{1}
Rf^{2}\,d\upsilon\\
\geq&\int_{M} (1-n\rho)\lambda_{1}
Rf^{2}\,d\upsilon-2a\int_{M}|\nabla f|^{2}\,d\upsilon\\
=&\int_{M} (1-n\rho)\lambda_{1} Rf^{2}\,d\upsilon-2a\lambda_{1}\\
=&\lambda_{1}\int_{M}f^{2}\{(1-n\rho)R-2a\}\,d\upsilon,
\endaligned\end{equation}
where we used the first assumption of Theorem 1.2.

 Since the eigenfunction of the first eigenvalue is not equal to $0$ almost everywhere along the Ricci-Bourguignon flow, by Lemma 2.3 and \eqref{1Th81},  we obtain
 \begin{equation}\label{1Th83}
\frac{d}{dt}\mu(t)\mid_{t=t_{0}}>0.
\end{equation}
Using the similar arguments as the proof of Theorem 1.2, we have $\lambda_{1}(t)$ is strictly monotone increasing in
$[0, T)$.
Therefore we finish the proof of Theorem 1.2.

As for the differentiability for $\lambda_{1}(t)$, since $\lambda_{1}(t)$ is increasing on the time interval
$[0, T )$ under curvature conditions of the theorem, by the classical Lebesgue's theorem (for
example, see Chap. 4 in \cite{Mukherjea1984}), it is easy to see that $\lambda_{1}(t)$ is differentiable almost everywhere
on $[0, T )$.

\noindent{\bf Remark 4.1.}  (1) In the course of proving Theorem
1.3, we do not use any differentiability of the first eigenvalue or
the corresponding eigenfunction of the Laplace operator under the
Ricci-Bourguignon flow.

(2) Using this method, we cannot get any monotonicity for higher order eigenvalues of the
Laplace operator under the Ricci-Bourguignon flow.

\section{Proof of Theorems 1.4 }

When $(M_{3}, g_{0})$ is a closed three manifold with positive Ricci curvature, the eigenvalues of the Laplacian along the
Ricci-Bourguignon flow diverges as $t \rightarrow T$. Now we finish
the proof of Theorem 1.3:

\noindent{\bf Proof of Theorems 1.4.} On a closed manifold
$M^{n}$, for any smooth functions $f$, using the celebrated Reilly
formula
$$\int_{M}|\nabla^{2}f|^{2}\,d\upsilon+
\int_{M}Ric(\nabla f, \nabla f)\,d\upsilon=\int_{M}(\Delta
f)^{2}\,d\upsilon$$ and Cauchy inequality
$$|\nabla^{2}f|^{2}\geq\frac{1}{n}(\Delta f)^{2},$$ we have the
inequality
$$\int_{M}Ric(\nabla f, \nabla
f)\,d\upsilon\leq\frac{n-1}{n}\int_{M}(\Delta f)^{2}\,d\upsilon.$$
When $\rho<\frac{1}{4}$, for any compact maximal solution of the Ricci-Bourguignon
flow on a closed three manifold with positive Ricci curvature there
exists $\varepsilon\leq\frac{1}{3}$ such that the condition
$$Ric\geq\varepsilon Rg$$ is preserved along the flow. For any time $t_{0}$, if $f$ is the eigenfunction of
the eigenvalue $\lambda$, then
$$\frac{2}{3}\lambda^{2}(t)\geq\int_{M}Ric(\nabla f, \nabla
f)\,d\upsilon\geq\varepsilon\int_{M}R|\nabla
f|^{2}\,d\upsilon\geq\varepsilon R_{min}(t)\lambda(t).$$ Since $t_{0} \in [0, T)$ is arbitrary, then
$$\lambda(t)\geq\frac{3}{2}\varepsilon R_{min}.$$
The thesis follows since $$\lim\limits_{t\rightarrow
T}R_{min}(t)=\infty,$$ see \cite{Fabrizio2007}.
Therefore we can finish the proof of Theorem 1.4.

\bibliographystyle{Plain}

\begin{thebibliography}{10}
\bibitem{Perelman2002} Perelman G. The entropy formula for the
Ricci flow and its geometric applications. ArXiv:0211159
\end{thebibliography}

FANQI ZENG\\
DEPARTMENT OF MATHEMATICS
TONGJI UNIVERSITY\\
SHANGHAI 200092
P.R. CHINA\\
E-mail: fanzeng10$@$126.com\\

QUN HE\\
DEPARTMENT OF MATHEMATICS
TONGJI UNIVERSITY\\
SHANGHAI 200092
P.R. CHINA\\
E-mail: hequn$@$tongji.edu.cn\\

BIN CHEN\\
DEPARTMENT OF MATHEMATICS
TONGJI UNIVERSITY\\
SHANGHAI 200092
P.R. CHINA\\
E-mail: Chenbin$@$tongji.edu.cn\\

\end{document}